\documentclass{article}
\usepackage{amsfonts}
\usepackage{makeidx}
\usepackage{amssymb}
\usepackage{amsmath}
\usepackage{amstext}
\usepackage{graphicx}

\newtheorem{theorem}{Theorem}

\newtheorem{corollary}[theorem]{Corollary}

\newtheorem{definition}[theorem]{Definition}

\newtheorem{proposition}[theorem]{Proposition}
\newtheorem{remark}[theorem]{Remark}

\input tcilatex
\begin{document}

\title{Tempered ultrafunctions}
\author{Vieri Benci\thanks{
Dipartimento di Matematica, Universit\`{a} degli Studi di Pisa, Via F.
Buonarroti 1/c, 56127 Pisa, ITALY and Centro Interdisciplinare "Beniamino
Segre", Accademia del Lincei. e-mail: \texttt{vieri.benci@unipi.it}}}
\maketitle

\begin{abstract}
Ultrafunctions are a particular class of functions defined on some
non-Archimedean field. They provide generalized solutions to functional
equations which do not have any solutions among the real functions or the
distributions. In this paper we introduce a new class of ultrafunctions,
called tempered ultrafunctions, which are somewhat related to the tempered
distributions and present some interesting peculiarities.

\medskip


\medskip

\end{abstract}

\tableofcontents

\bigskip

\section{Introduction}

In many circumstances, the notion of function is not sufficient to the needs
of a theory and it is necessary to generalize it. Generalized functions are
especially useful in making discontinuous functions more like smooth
functions. They are applied extensively, especially in physics and
engineering.

There is more than one recognized theory of \textit{generalized functions}.
We can recall, for example, the heuristic use of symbolic methods, called
operational calculus. A basic book on operational calculus was Oliver
Heaviside's Electromagnetic Theory of 1899 \cite{H99}. A very important
steps in this topic was the introduction of the weak derivative and of the
Dirac Delta function. The theory of Dirac and the theory of weak derivatives
where unified by Schwartz in the beautiful theory of distributions (see e.g. 
\cite{Sw51} and \cite{Sw66}), also thanks to the previous work of Leray and
Sobolev. Among people working in partial differential equations, the theory
of Schwartz has been accepted as definitive (at least until now), but other
notions of generalized functions have been introduced by Colombeau \cite%
{col85} and Mikio Sato \cite{sa59}, \cite{sa60}.

Having in mind the same purposes, in some recent papers, we have introduced
and studied the notion of ultrafunction (\cite{ultra}, \cite{belu2012}, \cite%
{belu2013}, \cite{milano}, \cite{algebra}, \cite{beyond}, \cite{topologia}, 
\cite{BLta}). Ultrafunctions are a particular class of functions defined on
a Non Archimedean field (we recall that a Non Archimedean field is an
ordered field which contains infinite and infinitesimal numbers).

This paper continues such studies and introduces a new class of
ultrafunctions called \textbf{tempered ultrafunctions}. This name comes from
the fact that tempered ultrafunctions are somewhat related to the tempered
distributions.

The theory of tempered ultrafunctions is based on four elements, $V_{\sigma
},(%
{{}^\circ}%
),\doint ,D:$

\begin{itemize}
\item $V_{\sigma }$ is a $\overline{\mathbb{R}}$-linear subspace of $%
\mathfrak{F}\left( \overline{\mathbb{R}},\overline{\mathbb{C}}\right) $
where $\overline{\mathbb{R}}$ is a Non-Archimedean field, $\overline{\mathbb{%
C}}=\overline{\mathbb{R}}\oplus i\overline{\mathbb{R}}$ and $\mathfrak{F}%
\left( \overline{\mathbb{R}},\overline{\mathbb{C}}\right) $ is the set of
all functions from $\overline{\mathbb{R}}\ $to $\overline{\mathbb{C}};$

\item 
\begin{equation*}
(%
{{}^\circ}%
):\mathcal{S}^{\prime }\left( \mathbb{R}\right) \rightarrow V_{\sigma }
\end{equation*}%
is an injective linear map defined on the set of tempered distributions $%
\mathcal{S}^{\prime }\left( \mathbb{R}\right) ;$

\item 
\begin{equation*}
\doint :V_{\sigma }\times V_{\sigma }\rightarrow \overline{\mathbb{C}}
\end{equation*}%
is an Hermitian bilinear form;

\item 
\begin{equation*}
D:V_{\sigma }\rightarrow V_{\sigma }
\end{equation*}%
is a $\overline{\mathbb{C}}$-linear operator which extends the
distributional derivative.
\end{itemize}

Some of the properties of the temperd ultrafunctions are described in the
following theorem:

\begin{theorem}
\label{one}The quadruplet $\left\{ V_{\sigma },(%
{{}^\circ}%
),\doint ,D\right\} $ satisfies the following properties:

\begin{enumerate}
\item \label{11}if $T\in \mathcal{S}^{\prime }\left( \mathbb{R}\right) ,$
then, $\forall \varphi \in \mathcal{S}\left( \mathbb{R}\right) $%
\begin{equation*}
\doint T%
{{}^\circ}%
(x)\varphi 
{{}^\circ}%
(x)dx=\left\langle T,\varphi \right\rangle ;
\end{equation*}

\item \label{22}if $f\in C^{0}\left( \mathbb{R}\right) $ and $\hat{f}\in
L^{1}\left( \mathbb{R}\right) ,$ then%
\begin{equation*}
\forall x\in \mathbb{R},\ \ f%
{{}^\circ}%
(x)\sim f(x);
\end{equation*}

\item \label{33}if $\varphi \in \mathcal{S}\left( \mathbb{R}\right) $ then%
\begin{equation*}
\forall x\in \mathbb{R},\ \ \varphi 
{{}^\circ}%
(x)=\varphi (x);
\end{equation*}

\item \label{44}if $T\in \mathcal{S}^{\prime }\left( \mathbb{R}\right) ,$
and $\partial $ is the distributional derivative, then%
\begin{equation*}
\left( \partial T\right) 
{{}^\circ}%
=DT%
{{}^\circ}%
\end{equation*}

\item \label{88}if $f\in C^{0}\left( \mathbb{R}\right) $ is a rapidly
decreasing function, then 
\begin{equation*}
\doint f%
{{}^\circ}%
(x)dx=\int f(x)dx
\end{equation*}

\item \label{55}if $u,v\in V_{\sigma },$%
\begin{equation*}
\doint Du(x)v(x)dx=-\doint u(x)Dv(x)dx;
\end{equation*}

\item \label{66}let $\mathcal{F}:\mathcal{S}^{\prime }\left( \mathbb{R}%
\right) \rightarrow \mathcal{S}^{\prime }\left( \mathbb{R}\right) $ denote
the distibutional Fourier transform; then%
\begin{equation*}
\mathcal{F}\left[ T\right] (k)=\frac{1}{\sqrt{2\pi }}\doint T%
{{}^\circ}%
(x)\left( e^{ikx}\right) 
{{}^\circ}%
dx
\end{equation*}
\end{enumerate}
\end{theorem}

However one of the most interesting peculiarities of the ultrafunctions
relies in the fact that they are functions and not functionals as the
disributions. This fact implies that also $\delta ^{2},$ $\sqrt{\delta }$
etc. are well defined functions in $\mathfrak{F}\left( \overline{\mathbb{R}},%
\overline{\mathbb{C}}\right) $.

\subsection{Notations and definitions\label{not}}

We use this section to fix some notations and to recall some definitions:

\begin{itemize}
\item $\mathfrak{F}\left( X,Y\right) $ denotes the set of all functions from 
$X$ to $Y;$

\item $\mathfrak{F}\left( E\right) =\mathfrak{F}\left( E,\mathbb{R}\right) ;$

\item $C^{k}\left( \mathbb{R}\right) $ denotes the set of functions in $%
C\left( \mathbb{R}\right) $ which have continuous derivatives up to the
order $k;$

\item $\mathcal{D}\left( \mathbb{R}\right) $ denotes the set of the
infinitely differentiable functions with compact support$;\ \mathcal{D}%
^{\prime }\left( \mathbb{R}\right) $ denotes the topological dual of $%
\mathcal{D}\left( \mathbb{R}\right) $, namely the set of distributions on $%
\mathbb{R};$

\item $\mathcal{S}\left( \mathbb{R}\right) $ denotes the set of the
infinitely differentiable functions rapidly decreasing with their
derivatives;$\ \mathcal{S}^{\prime }\left( \mathbb{R}\right) $ denotes the
topological dual of $\mathcal{S}\left( \mathbb{R}\right) $, namely the set
of tempered distributions on $\mathbb{R};$

\item $C_{\tau }^{0}\left( \mathbb{R}\right) $ denotes the set of continuous
functions slowly incresing (namely growing less than $x^{m}$ for some $m\in 
\mathbb{N}$);

\item if $\mathbb{K\ }$is a non-Archemedean field, then

\begin{itemize}
\item $\mathfrak{mon}(x)=\{y\in \mathbb{K}:x\sim y\};$

\item $\mathfrak{gal}(x)=\{y\in \mathbb{K}:x-y\ \ $is a finite number$\}.$
\end{itemize}
\end{itemize}

\section{Some notions of Non-Archimedean mathematics}

We believe that Non Archimedean Mathematics (NAM), namely, mathematics based
on Non Archimedean Fields is very interesting, very rich and, in many
circumstances, allows to construct models of the physical world in a more
elegant and simple way. In the years around 1900, NAM was investigated by
prominent mathematicians such as Du Bois-Reymond \cite{DBR}, Veronese \cite%
{veronese}, \cite{veronese2}, David Hilbert \cite{hilb}, and Tullio
Levi-Civita \cite{LC}, but then it has been forgotten until the '60s when
Abraham Robinson presented his Non Standard Analysis (NSA) \cite{rob}. We
refer to Ehrlich \cite{el06} for a historical analysis of these facts and to
Keisler \cite{keisler76} for a very clear exposition of NSA. Here, we will
construch a model of NAM based on \cite{benci99} and \cite{BDN2003}.

\subsection{Non Archimedean Fields\label{naf}}

In this section, we recall the basic definitions and facts regarding
non-Archimedean fields. In the following, ${\mathbb{K}}$ will denote an
ordered field. We recall that such a field contains (a copy of) the rational
numbers. Its elements will be called numbers.

\begin{definition}
Let $\mathbb{K}$ be an ordered field. Let $\xi \in \mathbb{K}$. We say that:

\begin{itemize}
\item $\xi $ is infinitesimal if, for all positive $n\in \mathbb{N}$, $|\xi
|<\frac{1}{n}$;

\item $\xi $ is finite if there exists $n\in \mathbb{N}$ such as $|\xi |<n$;

\item $\xi $ is infinite if, for all $n\in \mathbb{N}$, $|\xi |>n$
(equivalently, if $\xi $ is not finite).
\end{itemize}
\end{definition}

\begin{definition}
An ordered field $\mathbb{K}$ is called Non-Archimedean if it contains an
infinitesimal $\xi \neq 0$.
\end{definition}

It's easily seen that all infinitesimal are finite, that the inverse of an
infinite number is a nonzero infinitesimal number, and that the inverse of a
nonzero infinitesimal number is infinite.

\begin{definition}
A superreal field is an ordered field $\mathbb{K}$ that properly extends $%
\mathbb{R}$.
\end{definition}

It is easy to show, due to the completeness of $\mathbb{R}$, that there are
nonzero infinitesimal numbers and infinite numbers in any superreal field.
Infinitesimal numbers can be used to formalize a new notion of "closeness":

\begin{definition}
\label{def infinite closeness} We say that two numbers $\xi, \zeta \in {%
\mathbb{K}}$ are infinitely close if $\xi -\zeta $ is infinitesimal. In this
case, we write $\xi \sim \zeta $.
\end{definition}

Clearly, the relation "$\sim $" of infinite closeness is an equivalence
relation.

\begin{theorem}
If $\mathbb{K}$ is a superreal field, every finite number $\xi \in \mathbb{K}
$ is infinitely close to a unique real number $r\sim \xi $, called the 
\textbf{shadow} or the \textbf{standard part} of $\xi $.
\end{theorem}

Given a finite number $\xi $, we denote its shadow as $sh(\xi )$, and we put 
$sh(\xi )=+\infty $ ($sh(\xi )=-\infty $) if $\xi \in \mathbb{K}$ is a
positive (negative) infinite number.\newline

\begin{definition}
Let $\mathbb{K}$ be a superreal field, and $\xi \in \mathbb{K}$ a number.
The \label{def monad} monad of $\xi $ is the set of all numbers that are
infinitely close to it:%
\begin{equation*}
\mathfrak{m}\mathfrak{o}\mathfrak{n}(\xi )=\{\zeta \in \mathbb{K}:\xi \sim
\zeta \},
\end{equation*}%
and the galaxy of $\xi $ is the set of all numbers that are finitely close
to it: 
\begin{equation*}
\mathfrak{gal}(\xi )=\{\zeta \in \mathbb{K}:\xi -\zeta \ \text{is\ finite}\}
\end{equation*}
\end{definition}

By definition, it follows that the set of infinitesimal numbers is $%
\mathfrak{mon}(0)$ and that the set of finite numbers is $\mathfrak{gal}(0)$.

\bigskip

\subsection{Sigma-convergence}

When we take a limit of a sequence $\varphi (n)$ for $n\rightarrow \infty ,$
the family of neighborhoods of $\infty $ is the Frechet filter, namely the
family of cofinite sets. In order to realize our program we need a finer
topology based on an ultrafilter.

\begin{definition}
\label{pula}An ultrafilter $\sigma $ on $\mathbb{N}$ is a filter which
satisfies the following property:\ if $A\cup B=\mathbb{N}$, then 
\begin{equation}
A\in \mathfrak{\sigma \ }\text{or }B\in \sigma ;  \label{quaqua}
\end{equation}

We will refer to the sets in $\sigma $ as \textbf{qualified sets}. If a
property $P(n)$ holds for every $n$ in a qualified set we will say that $P$
holds a.e. (almost everywhere).
\end{definition}

Now we define a topology on $\mathbb{N}\cup \left\{ \sigma \right\} $ as
follows: $N$ is a neighborhood of $\sigma $ if and only if 
\begin{equation}
N=Q\cup \left\{ \sigma \right\} \ \ \ \ \text{with}\ \ \ \ Q\in \sigma
\label{alfa}
\end{equation}

Thus $\mathbb{N}\cup \left\{ \sigma \right\} $ and $\mathbb{N}\cup \left\{
\infty \right\} \mathbb{\ }$are similar spaces but have different
topologies: in fact $\mathbb{N}\cup \left\{ \sigma \right\} $ has a finer
topology than $\mathbb{N}\cup \left\{ \infty \right\} $. In any case you may
think of $\sigma $ as a "point at infinity".

\bigskip

\begin{definition}
Given a topological space $X,$ we say that a sequence $\varphi :\mathbb{N}%
\rightarrow X$ is $\sigma $-\textbf{convergent }to $L\in X$ if and only if
for any neighborhood $N$ of $L,\ \exists Q\in \sigma $\ such that, 
\begin{equation}
\forall n\in Q,\ \varphi _{n}\in N.  \label{limite}
\end{equation}%
We will write%
\begin{equation*}
\underset{n\uparrow \sigma }{\lim }\ \varphi _{n}=L
\end{equation*}
\end{definition}

Sometimes, the $\sigma $-convergence is called convergence along an
ultrafilter. Notice that this topology satisfies this interesting property:

\begin{proposition}
\label{nino}If the sequence $\varphi _{n}\in X$ has a converging
subsequence, then it is $\sigma $-\textbf{converging}.
\end{proposition}

\textbf{Proof}: Suppose that the net $\varphi _{n}$ has a converging subnet
to $L\in \mathbb{R}$. We a neighborhood $N$ of $L$ and we have to prove that 
$Q_{N}\in \sigma $ where%
\begin{equation*}
Q_{N}=\left\{ \lambda \in \mathbb{N}\ |\ \varphi _{n}\in N\right\} .
\end{equation*}%
We argue indirectly and we assume that 
\begin{equation*}
Q_{N}\notin \sigma
\end{equation*}%
Then, by (\ref{quaqua}), $N=\mathbb{N}\backslash \left( Q_{N}\cap X\right)
\in \sigma $ and hence%
\begin{equation*}
\forall n\in N,\ \varphi _{n}\notin N.
\end{equation*}%
This contradict the fact that $\varphi _{n}$ has a subsequence which
converges to $L.$

$\square $

\subsection{The fields of tempered hyperreal and hypercomplex numbers}

In this section, we will construct a new Archimedean field, strictly related
to $\mathbb{R}^{\mathbb{\ast }}$ which we will call field of the \textbf{%
tempered hyperreal numbers }and we will denote by $\overline{\mathbb{R}}$.
We recall that $\mathbb{R}^{\mathbb{\ast }},$ the field of hyperreal
numbers, is the basic field in Nonstandard Analysis (see e.g. \cite{rob} or 
\cite{keisler76} for a very clear exposition; see \cite{benci95}, \cite%
{benci99}, or \cite{BDN2003} for an exposition closer to the content of this
paper). However, even if some ideas on Nonstandard Analysis are used, it is
not necessary for the reader to be familiar with it in order to read this
paper.

\bigskip 

We define the family of \textbf{slowly increasing sequences} as follows 
\begin{equation*}
\mathfrak{F}_{\tau }\left( \mathbb{N},\mathbb{R}\right) =\left\{ \varphi \in 
\mathfrak{F}\left( \mathbb{N},\mathbb{R}\right) \ |\ \ \exists \ k\in 
\mathbb{N},\ |\varphi _{n}|\leq k+n^{k}\right\} \ 
\end{equation*}%
We say that a sequence $\varphi _{n}$ is \textbf{rapidly increasing }if,%
\textbf{\ }$\forall k\in \mathbb{N},$ $\varphi _{n}\geq n^{k}$ for every $%
n\geq \overline{n}\in \mathbb{N}$.

Moreover the set of \textbf{rapidly decreasing sequences }(with respect to $%
\sigma $) is defined as follows:%
\begin{equation*}
\mathfrak{I}\left( \mathbb{N},\mathbb{R}\right) =\left\{ \varphi \in 
\mathfrak{F}\left( \mathbb{N},\mathbb{R}\right) \ |\ \ \forall k\in \mathbb{N%
},\exists Q\in \sigma ,\forall n\in Q,\ |\varphi _{n}|\leq n^{-k}\right\} \ 
\end{equation*}%
More in general, given an unbounded set $E\subset \mathbb{R}$, we define the
family of \textbf{slowly increasing functions} as follows: 
\begin{equation*}
\mathfrak{F}_{\tau }\left( E,\mathbb{R}\right) =\left\{ \varphi \in 
\mathfrak{F}\left( E,\mathbb{R}\right) \ |\ \exists \ k\in \mathbb{N},\
\forall x\in E,\ |\varphi (x)|\leq k+|x|^{k}\right\} 
\end{equation*}%
and similarly we can define the rapidly increasing functions and rapidly
decreasing functions.

We have the following result:

\begin{theorem}
\label{nuovo}\textit{There exists an ordered field} $\overline{\mathbb{R}}%
\supset \mathbb{R}$\textit{\ (}$\overline{\mathbb{R}}\neq \mathbb{R}$\textit{%
) such that }

\begin{enumerate}
\item \textit{Every sequence }$\varphi \in \mathfrak{F}_{\tau }\left( 
\mathbb{N},\mathbb{R}\right) $\textit{\ has a unique limit } 
\begin{equation*}
\underset{n\uparrow \sigma }{\lim }\ \varphi _{n}=L\in \overline{\mathbb{R}}.
\end{equation*}%
which will be called $\sigma $-limit. \textit{Moreover every}\emph{\ }$\xi
\in \overline{\mathbb{R}}$\textit{\ is the }$\sigma $-\textit{limit\ of some
sequence }$\varphi :\mathbb{N}\rightarrow \mathbb{R}$\emph{.}

\item If $\varphi \in \mathfrak{I}\left( \mathbb{N},\mathbb{R}\right) $,
then 
\begin{equation*}
\underset{n\uparrow \sigma }{\lim }\ \varphi _{n}=0
\end{equation*}

\item \textit{For all }$\varphi ,\psi :\mathbb{N}\rightarrow \mathbb{R}$%
\emph{: }%
\begin{eqnarray*}
\underset{n\uparrow \sigma }{\lim }\ \varphi _{n}+\ \underset{n\uparrow
\sigma }{\lim }\ \psi _{n} &=&\ \underset{n\uparrow \sigma }{\lim }\left(
\varphi _{n}+\psi _{n}\right) ; \\
\underset{n\uparrow \sigma }{\lim }\ \varphi _{n}\cdot \ \underset{n\uparrow
\sigma }{\lim }\ \psi _{n} &=&\ \underset{n\uparrow \sigma }{\lim }\left(
\varphi _{n}\cdot \psi _{n}\right) .
\end{eqnarray*}
\end{enumerate}
\end{theorem}

\textbf{Proof. }Since $\mathfrak{I}\left( \mathbb{N},\mathbb{R}\right) $ is
an ideal in $\mathfrak{F}_{\tau }\left( \mathbb{N},\mathbb{R}\right) ,\ $we
have that%
\begin{equation*}
\overline{\mathbb{R}}=\frac{\mathfrak{F}_{\tau }\left( \mathbb{N},\mathbb{R}%
\right) }{\mathfrak{I}\left( \mathbb{N},\mathbb{R}\right) }
\end{equation*}%
is a ring and it is not difficult to prove that it is an ordered field
(namely $\mathfrak{I}\left( \mathbb{N},\mathbb{R}\right) $ is a maximal
ideal). If we set 
\begin{equation*}
\underset{n\uparrow \sigma }{\lim }\ \varphi _{n}=\left[ \varphi \right] 
\end{equation*}%
we have that every sequence has a unique limit and if we identify a real
number $c\in \mathbb{R}$ with the equivalence class of the constant sequence 
$\left[ c\right] ,$ then $\mathbb{R}\subset \overline{\mathbb{R}}.$

$\square $

\bigskip We have the following result:

\begin{theorem}
\label{nuovo1}There is a topology on $\overline{\mathbb{R}}$ consistent with
the $\sigma $-\textit{limit. This topology will be called }$\sigma $%
-topology.
\end{theorem}

\textbf{Proof}. Probably, the simplest way to define this topology is giving
the closure operator as follows: given a set $E\subset \overline{\mathbb{R}}%
, $ we set,%
\begin{eqnarray*}
D\left( E\right) &=&\left\{ \xi \in \overline{\mathbb{R}}\backslash \mathbb{R%
}\ |\ \xi =\ \underset{n\uparrow \sigma }{\lim }\ \varphi _{n}\ \text{where}%
\ \forall n\in \mathbb{N},\ \varphi _{n}\in E\cap \mathbb{R}\right\} \ \ if\
\ E\cap \mathbb{R}\neq \varnothing \\
D(E) &=&\varnothing \ \ \ if\ \ E\subset \overline{\mathbb{R}}\backslash 
\mathbb{R}=\varnothing \\
\overline{E} &=&E\cup D\left( E\right)
\end{eqnarray*}

We will show that this operator satisfies the Kuratowski closure axioms:

\begin{enumerate}
\item Preservation of the empty set: $\overline{\varnothing }=\varnothing .$

\item Extensivity: $E\subset \overline{E}$

\item Preservation of Binary Union: $\overline{E\cup F}=\overline{E}\cup 
\overline{F}$

\item Idempotence: $\overline{\overline{E}}=\overline{E}$
\end{enumerate}

1 and 2 are trivial. Let us prove $3$.\ By Def. \ref{pula}, if $\varphi
_{n}\in \left( E\cup F\right) \cap \mathbb{R},$ there exists a qualified set 
$Q_{\varphi }$ such that $\forall n\in Q_{\varphi },\ \varphi _{n}\in E\cap 
\mathbb{R}$ or $\varphi _{n}\in F\cap \mathbb{R};$ then%
\begin{eqnarray*}
D(E\cup F) &=&\left\{ \underset{n\uparrow \sigma }{\lim }\ \varphi _{n}\ |\
\forall n\in \mathbb{N},\ \varphi _{n}\in \left( E\cup F\right) \cap
X\right\} \\
&=&\left\{ \underset{n\uparrow \sigma }{\lim }\ \varphi _{n}\ |\ \forall
n\in Q_{\varphi },\ \varphi _{n}\in E\cap X\right\} \cup \left\{ \underset{%
n\uparrow \sigma }{\lim }\ \varphi _{n}\ |\ \forall n\in Q_{\varphi },\
\varphi _{n}\in F\cap X\right\} \\
&=&D(E)\cup D(F)
\end{eqnarray*}

So, we have that%
\begin{eqnarray*}
\overline{E\cup F} &=&E\cup F\cup D(E\cup F) \\
&=&E\cup F\cup D(E)\cup D(F) \\
&=&\overline{E}\cup \overline{F}
\end{eqnarray*}

Now let us prove $4$: Since $D\left( E\right) \subset \overline{\mathbb{R}}%
\backslash \mathbb{R},$ we have that $D(D\left( E\right) )=\varnothing ,$
and hence $\overline{D\left( E\right) }=D\left( E\right) ;$ thus 
\begin{equation*}
\overline{\overline{E}}=\overline{E\cup D\left( E\right) }=\overline{E}\cup 
\overline{D\left( E\right) }=\left[ E\cup D\left( E\right) \right] \cup
D\left( E\right) =\overline{E}.
\end{equation*}

$\square $

From now on, unless differently stated, the notation 
\begin{equation*}
\underset{n\uparrow \sigma }{\lim }\ \varphi _{n}
\end{equation*}%
will denote the $\sigma $-limit of the sequence $\varphi _{n}$ where the
target space is $\overline{\mathbb{R}}$ with its topology; the notation%
\begin{equation*}
\underset{n\rightarrow \infty }{\lim }\ \varphi _{n}
\end{equation*}%
will denote the usual limit, where the target space is $\mathbb{R}$ with the
usual topology.

\begin{theorem}
\label{beatrice}Given a sequence $\varphi _{n}\in \mathbb{R}$, then we have
one of the mutually exclusive possibility:

\begin{enumerate}
\item there is $Q\in \sigma $ such that $\varphi _{n}|_{Q}$ is slowly
increasing; in this case we have that%
\begin{equation*}
\underset{n\uparrow \sigma }{\lim }\ \varphi _{n}\in \overline{\mathbb{R}}
\end{equation*}

\item there is $Q\in \sigma $ such that $\varphi _{n}|_{Q}$ is rapidly
increasing; in this case we will write%
\begin{equation*}
\underset{n\uparrow \sigma }{\lim }\ \varphi _{n}=+\overline{\infty }
\end{equation*}

\item there is $Q\in \sigma $ such that $-\varphi _{n}|_{Q}$ is rapidly
increasing; in this case we will write%
\begin{equation*}
\underset{n\uparrow \sigma }{\lim }\ \varphi _{n}=-\overline{\infty }
\end{equation*}
\end{enumerate}
\end{theorem}

\textbf{Proof: }It is an immediate consequence of the maximality of $\sigma
. $

$\square $\bigskip

Using the above notation we have that $\left\{ \overline{\mathbb{R}},+%
\overline{\infty },-\overline{\infty }\right\} $ is the compactification of $%
\overline{\mathbb{R}}$ analogous of the extended real line $\left\{ \mathbb{R%
},+\infty ,-\infty \right\} .$

The following theorem whose proof is left to the reader shows some relation
existing between the usual limit and the $\sigma $-limit.

\begin{theorem}
\label{peso}If $\varphi =\ \underset{n\rightarrow \infty }{\lim }\ \varphi
_{n},$ then%
\begin{equation*}
\underset{n\uparrow \sigma }{\lim }\ \varphi _{n}\sim \varphi .
\end{equation*}%
Moreover, if $\varphi -\varphi _{n}$ is a rapidly decreasing sequence, then%
\begin{equation*}
\underset{n\uparrow \sigma }{\lim }\ \varphi _{n}=\varphi .
\end{equation*}
\end{theorem}

Since in the next section, we will deal with complex valued function, the
following trivial definition is needed.

\begin{definition}
The complexification of $\overline{\mathbb{R}}$ is called field of the
tempered hypercomplex numbers and it is denoted by $\overline{\mathbb{C}}$
\end{definition}

\bigskip

\section{The space of tempered ultrafunctions $V_{\protect\sigma }$}

\bigskip

\subsection{Definition of tempered ultrafunctions}

We set%
\begin{equation}
\beta _{n}:=n\sqrt{\pi };\ \eta _{n}:=\frac{\sqrt{\pi }}{n}  \label{gina}
\end{equation}%
\begin{equation}
\Sigma _{n}=\left\{ l\eta _{n}\ |\
l=-n^{2},-n^{2}+1,-n^{2}+2,....,n^{2}-2,n^{2}-1\right\}  \label{rina}
\end{equation}%
\begin{equation}
V_{n}=\text{Span}\left\{ e^{ikx}\ |\ k\in \Sigma \right\} .  \label{pina}
\end{equation}

So, $V_{n}$ satisfies the following properties:

\begin{itemize}
\item $V_{n}$ is a vector space of dimension $2n^{2}$

\item if $f\in V_{n},$ $f$ is periodic of period $2\beta _{n}$
\end{itemize}

We now set%
\begin{equation*}
\Sigma =\left\{ \xi \in \overline{\mathbb{R}}\ |\ \exists \varphi \in 
\mathfrak{F}_{\tau }\left( \mathbb{N},\mathbb{R}\right) ,\ \varphi _{n}\in
\Sigma _{n},\ \xi =\ \underset{n\uparrow \sigma }{\lim }\ \varphi
_{n}\right\}
\end{equation*}

\bigskip

\begin{definition}
Let $E\subset \mathbb{R}$ and let $\overline{E}$ denote its closure with
respect to the $\sigma $-topology on $\overline{\mathbb{R}}$. A function%
\begin{equation*}
u:\overline{E}\rightarrow \overline{\mathbb{R}}
\end{equation*}%
is called $\sigma $-limit function if there is a sequence of functions $%
u_{n}\in \mathfrak{F}\left( E,\mathbb{R}\right) $ such that $\forall \xi \in 
\overline{E},$%
\begin{equation}
u(\xi )=\mathfrak{\ }\underset{n\uparrow \sigma }{\lim }\ u_{n}(x_{n})\ 
\label{lim}
\end{equation}%
where%
\begin{equation*}
\underset{n\uparrow \sigma }{\lim }\ x_{n}=\xi ;\ x_{n}\in E.
\end{equation*}%
The set of $\sigma $-limit functions will be denoted by $\mathfrak{S}\left( 
\overline{E},\overline{\mathbb{R}}\right) $. The set of tempered
ultrafunction $V_{\sigma }\subset \mathfrak{F}\left( \overline{\mathbb{R}},%
\overline{\mathbb{R}}\right) $ is the subset of the $\sigma $-limit
functions such that $\forall n\in \mathbb{N}$, 
\begin{equation}
u_{n}\in V_{n}.  \label{extra}
\end{equation}
\end{definition}

\bigskip

In the following, in order to denote the limit (\ref{lim}), we will use the
following shorthand notation:%
\begin{equation*}
u=\mathfrak{\ }\underset{n\uparrow \sigma }{\lim }\ u_{n}.
\end{equation*}

\begin{remark}
If $f_{n}$ is a sequence of real functions, and $\xi =\ \underset{n\uparrow
\sigma }{\lim }\ x_{n},$ then, by Th. \ref{beatrice}, 
\begin{equation}
f_{\sigma }(\xi ):=\mathfrak{\ }\underset{n\uparrow \sigma }{\lim }\
f_{n}(x_{n})\ 
\end{equation}%
is always defined. However, if the sequence $|f_{n}(x_{n})|$ grows rapidly,
we have that $f_{\sigma }(\xi )=+\overline{\infty }\ $or$\ f_{\sigma }(\xi
)=-\overline{\infty }.$
\end{remark}

Let us see an example of ultrafunction:

\bigskip

\textbf{The exponential ultrafunction. }If $x\in \overline{\mathbb{R}},\
k\in \Sigma $, we set%
\begin{equation*}
e^{ikx}=\ \underset{n\uparrow \sigma }{\lim }\ e^{ik_{n}x_{n}}
\end{equation*}%
where $\underset{n\uparrow \sigma }{\lim }\ x_{n}=x$ and $\underset{%
n\uparrow \sigma }{\lim }\ k_{n}=k$. If $x,k\in \mathbb{R}$, then $e^{ikx}$
assume the same real values than the function $e^{ikx}$; for this reason,
they are denoted by the same symbols. In particular, for $k=0$, we will
denote by $1$ the tempered ultrafunction identically equal to (the tempered
hyperreal number) 1.

\subsection{Basic operation with the tempered ultrafunctions}

\bigskip

\subsubsection{Derivative}

The derivative of a tempered ultrafunction $u=\mathfrak{\ }\underset{%
n\uparrow \sigma }{\lim }\ u_{n},\ u_{n}\in V_{n}$ is defined as follows:%
\begin{equation*}
Du=\mathfrak{\ }\underset{n\uparrow \sigma }{\lim }\ u_{n}^{\prime }.
\end{equation*}

It is immediate to verify that the derivative satisfies the basic properties:

\bigskip

\begin{theorem}
The derivative satisfies the following properties:

\begin{enumerate}
\item \textbf{Linearity:}$\ D:V_{\sigma }\rightarrow V_{\sigma }$ is a
linear operator over the field $\overline{\mathbb{C}},$

\item \textbf{Leibnitz Rule:} if $u,v,D\left( uv\right) ,Duv,uDv\in
V_{\sigma }$ then,%
\begin{equation*}
D\left( uv\right) =Duv+uDv
\end{equation*}
\end{enumerate}
\end{theorem}

\bigskip

\begin{remark}
Clearly, it is possible to define the derivative also of a $\sigma $-limit
functions $u=\mathfrak{\ }\underset{n\uparrow \sigma }{\lim }\ u_{n}$
provides that each $u_{n}$ is differentiable.
\end{remark}

\bigskip

\subsubsection{Inner product}

\bigskip

The inner product between two $\sigma $-limit functions is defined as
follows:%
\begin{equation*}
(u|v)=\ \underset{n\uparrow \sigma }{\lim }\int_{-\beta _{n}}^{\beta
_{n}}u_{n}(x)\overline{v_{n}(x)}dx
\end{equation*}%
So $V_{\sigma }$ is an Euclidean space over the field of the tempered
hypercomplex numbers $\overline{\mathbb{C}}.$ The inner product induces on
the algebra of ultrafunctions the following "Euclidean" norm%
\begin{equation*}
\left\Vert u\right\Vert =\sqrt{(u|u)}
\end{equation*}%
which takes values in $\overline{\mathbb{R}}$.

\bigskip

\subsubsection{Integration}

\bigskip

The integral of a $\sigma $-limit function%
\begin{equation*}
u=\mathfrak{\ }\underset{n\uparrow \sigma }{\lim }\ u_{n},\ u_{n}\in
L_{loc}^{1}
\end{equation*}%
is defined as follows:%
\begin{equation*}
\doint u(x)dx=\ \underset{n\uparrow \sigma }{\lim }\int_{-\beta _{n}}^{\beta
_{n}}u_{n}(x)dx.
\end{equation*}%
Notice that if $n\mapsto \int_{-\beta _{n}}^{\beta _{n}}u_{n}(x)dx$ is a
slowly increasing function, the above limit exists in $\overline{\mathbb{R}}$
and we say that the integral converges; otherwise we say that the integral
divergis to $\overline{+\infty }$ or to $\overline{-\infty }$.

Using this definition, we have that%
\begin{equation}
(u|v)=\doint u_{n}(x)\overline{v_{n}(x)}dx  \label{pippa}
\end{equation}%
where $\overline{z}$ denotes the complex conjugate of $z.$

\bigskip

It is immediate to prove the following theorem:

\begin{theorem}
\label{pizzi}The integral satisfies the following properties:

\begin{enumerate}
\item $\doint :V_{\sigma }\rightarrow \overline{\mathbb{C}}$ is a linear
functional over the field $\overline{\mathbb{C}}$;

\item if $u,v\in V_{\sigma },\ \doint Du(x)v(x)dx=-\doint u(x)Dv(x)dx$

\item if $u\in V_{\sigma },\ \doint Du(x)dx=0$
\end{enumerate}
\end{theorem}

\subsubsection{Hyperfinite sums}

\begin{definition}
A set $\Gamma \subset \overline{\mathbb{R}}$ is called hyperfinite is there
is a sequence of finite sets $\Gamma _{n}\subset \mathbb{R}$ such that%
\begin{equation*}
\Gamma =\left\{ \xi \in \overline{\mathbb{R}}\ |\ \xi =\ \underset{n\uparrow
\sigma }{\lim }\ \varphi _{n},\ \text{and}\ \forall n\in \mathbb{N},\
\varphi _{n}\in \Gamma _{n}\right\}
\end{equation*}
\end{definition}

\begin{definition}
If $\Gamma $ is a hyperfinite set and%
\begin{equation*}
u:\Gamma \rightarrow \overline{\mathbb{R}};\ u=\ \underset{n\uparrow \sigma }%
{\lim }\ u_{n};\ u_{n}:\Gamma _{n}\rightarrow \overline{\mathbb{R}},
\end{equation*}%
is a $\sigma $-limit function, we set 
\begin{equation*}
\sum_{k\in \Gamma }u\left( k\right) =\ \underset{n\uparrow \sigma }{\lim }\
\sum_{k\in \Gamma _{n}}u_{n}\left( k\right)
\end{equation*}%
\bigskip
\end{definition}

\bigskip

$\sum_{k\in \Gamma }u\left( k\right) $ will be called the hyperfinite sum of
the $u(k)^{\prime }s.$

More in general, a family of $\sigma $-limit functions $\left\{
u_{k}\right\} _{k\in \Gamma }$ will be called hyperfinite if $\Gamma \subset 
\overline{\mathbb{R}}$ is a hyperfinite set and the functional%
\begin{equation*}
k\mapsto u_{k}
\end{equation*}%
is a $\sigma $-limit function. In this case, if $x=\mathfrak{\ }\underset{%
n\uparrow \sigma }{\lim }\ x_{n},$ we set%
\begin{equation*}
u(x)=\ \underset{n\uparrow \sigma }{\lim }\ \sum_{k\in \Gamma
_{n}}u_{k}(x_{n})
\end{equation*}%
It easy to check, just taking the appropriate limits from the finite case,
that the hyperfinite sum commute with the integral, namely 
\begin{equation*}
\sum_{j\in \Gamma }\left( \int u_{j}(x)dx\right) =\int \left( \sum_{j\in
\Gamma }u_{j}(x)\right) dx.
\end{equation*}

\section{Bases}

\bigskip

\subsection{The trigonometric basis}

\bigskip

\begin{definition}
An hyperfinite set of ultrafunctions $\left\{ e_{j}(x)\right\} _{j\in \Gamma
}$ is called a basis for $V_{\sigma }$ if every ultrafunction $u\in
V_{\sigma }$ can be written in a unique way by mean of an hyperfinite sum:%
\begin{equation*}
u(x)=\sum_{j\in \Gamma }u_{j}e_{j}(x)
\end{equation*}%
A bases $\left\{ e_{j}(x)\right\} _{j\in \Gamma }$ is said orthonomal if%
\begin{equation*}
\doint e_{j}(x)e_{k}(x)dx=\delta _{k}^{j}.
\end{equation*}
\end{definition}

\bigskip

By our construction we have that%
\begin{equation*}
\left\{ \frac{e^{ikx}}{\sqrt{2\beta _{n}}}\right\} _{k\in \Sigma _{n}}
\end{equation*}%
is an orthonormal basis for $V_{n}.$ Hence it is immediate to see that%
\begin{equation}
\left\{ \frac{e^{ikx}}{\sqrt{2\beta }}\right\} _{k\in \Sigma };\ \ \beta =\ 
\underset{n\uparrow \sigma }{\lim }\ \beta _{n}  \label{pilla}
\end{equation}%
is an othonormal basis for $V_{\sigma }$. We will refer to it as the \textbf{%
trigonometric basis}. So we have the following result:

\begin{theorem}
Any ultrafunction $u\in V_{\sigma }$ can be represented as follows:%
\begin{equation}
u(x)=\frac{1}{2\beta }\sum_{k\in \Sigma }\left( \doint u(t)e^{-ikt}dt\right)
e^{ikx}=\frac{\eta }{2\pi }\sum_{k\in \Sigma }\left( \doint
u(t)e^{-ikt}dt\right) e^{ikx}  \label{milva}
\end{equation}
\end{theorem}

Notice that the functions $e^{ikx}$ are the eigenfunctions of the operator $%
-iD_{x}$ with eigevalue $k$; this operator restricted to $V_{\sigma }$ is an
Hermitian operator and $\left\{ \left( 2\beta \right) ^{-1/2}e^{ikx}\right\}
_{k\in \Sigma _{n}}$ is the relative orthonormal basis of eigenvalues.

Thus the above formula can be seen as the ultrafunctions variant of the
relative spectral formula for selfadjoint operators. Also, as we will see
later, it is strictly related to the Fourier transform.

\subsection{The delta ultrafunction}

\bigskip

\begin{definition}
\label{dd}Given a number $q\in \overline{\mathbb{R}},$ we denote by $\delta
_{q}(x)$ an ultrafunction in $V_{\sigma }$ such that 
\begin{equation}
\forall v\in V_{\sigma },\ \doint v(x)\delta _{q}(x)dx=v(q)
\label{deltafunction}
\end{equation}%
$\delta _{q}(x)$ is called Delta (or the Dirac) ultrafunction concentrated
in $q$.
\end{definition}

Let us see the first properties of the Delta ultrafunctions:

\begin{theorem}
\label{delta} We have the following properties:

\begin{enumerate}
\item \label{dl1} for every $q\in \overline{\mathbb{R}}$ there exists an
unique Delta ultrafunction concentrated in $q;$

\item $\overline{\delta _{q}(x)}=\delta _{q}(x)$

\item for every $a,\ b\in \overline{\mathbb{R}},\ \delta _{a}(b)=\delta
_{b}(a);$

\item $\left\Vert \delta _{q}\right\Vert ^{2}=\delta _{q}(q).$
\end{enumerate}
\end{theorem}

\textbf{Proof.} (1) Let $\left\{ e_{j}\right\} _{j\in \Gamma }$ be an
orthonormal basis of $V_{\sigma }.$ We set 
\begin{equation}
\delta _{q}(x)=\sum_{j\in \Gamma }e_{j}(q)\overline{e_{j}(x)}\ \ \ \text{%
for\ any}\ q\in \mathbb{R}  \label{bimba}
\end{equation}

Let us prove that $\delta _{q}(x)$ actually satisfies (\ref{deltafunction}).
Let $v(x)=\sum_{j\in \Gamma }v_{j}e_{j}(x)$ be any ultrafunction and $q\in 
\mathbb{R}$. Then%
\begin{eqnarray*}
\doint v(x)\delta _{q}(x)dx &=&\int \left( \sum_{j\in \Gamma
}v_{j}e_{j}(x)\right) \left( \sum_{k\in \Gamma }e_{k}(q)\overline{e_{k}(x)}%
\right) dx= \\
&=&\sum_{j,k\in \Gamma }v_{j}e_{k}(q)\doint e_{j}(x)\overline{e_{k}(x)}dx= \\
&=&\sum_{j,k\in \Gamma }v_{j}e_{k}(q)\delta _{j}^{q}=\sum_{k\in \Gamma
}v_{k}e_{k}(q)=v(q).
\end{eqnarray*}

So $\delta _{q}(x)$, defined by (\ref{bimba}), is a Delta ultrafunction
centered in $q$.

It is unique: if $h_{q}(x)$ is another Delta ultrafunction centered in $q$
then for every $x\in \overline{\mathbb{R}},$ we have:%
\begin{eqnarray*}
\delta _{q}(x)-h_{q}(x) &=&\doint (\delta _{q}(y)-h_{q}(y))\delta
_{x}(y)dy=\doint \delta _{x}(y)\delta _{q}(y)dy-\doint \delta
_{x}(y)h_{q}(y)dy \\
&=&\delta _{x}(q)-\delta _{x}(q)=0
\end{eqnarray*}%
and hence $\delta _{q}(x)=h_{q}(x)$ for every $x\in \overline{\mathbb{R}}.$

(2) We may assume that the $e_{j}(q)$ are real functions in the sense that $%
\overline{e_{j}(q)}=e_{j}(q).$ Since, the $e_{j}(q)$ are real functions, we
have that $\overline{\delta _{q}(x)}=\delta _{q}(x).$

(3)$\ \delta _{a}\left( b\right) =\doint \delta _{a}(x)\delta _{b}(x)\
dx=\delta _{b}\left( a\right) .$

(4) $\left\Vert \delta _{q}\right\Vert ^{2}=\doint \delta _{q}(x)\delta
_{q}(x)=\delta _{q}(q)$.\newline

$\square $

\bigskip

\bigskip

\subsection{The canonical basis}

\bigskip

\bigskip

\begin{theorem}
\label{deltadelta}For any $a,b\in \Sigma $ 
\begin{equation*}
\doint \delta _{a}(x)\delta _{b}(x)dx=\frac{\delta _{b}^{a}}{\eta }
\end{equation*}%
where $\delta _{b}^{a}$ denotes the $\delta $ of Kronecker and%
\begin{equation*}
\eta =\ \underset{n\uparrow \sigma }{\lim }\eta _{n}
\end{equation*}
\end{theorem}

\textbf{Proof.} We have%
\begin{equation*}
\doint \delta _{a}(x)\delta _{b}(x)dx=\ \underset{n\uparrow \sigma }{\lim }%
\int_{-\beta _{n}}^{\beta _{n}}\delta _{n,a_{n}}(x)\delta _{n,b_{n}}(x)dx
\end{equation*}%
where $\delta _{n,a_{n}}\ $and$\ \delta _{n,b_{n}}$ are the approximations
of $\delta _{a}\ $and $\delta _{b}$ defined, using the orthonormal basis 
\begin{equation*}
\left\{ \frac{e^{ikx}}{\sqrt{2\beta _{n}}}\right\} _{k\in \Sigma _{n}},
\end{equation*}%
as follows%
\begin{equation}
\delta _{n,q}(x)=\frac{1}{2\beta _{n}}\sum_{k\in \Sigma _{n}}e^{ikq}e^{-ikx}.
\label{pirlo}
\end{equation}

Then we have 
\begin{eqnarray*}
\doint \delta _{a}(x)\delta _{b}(x)dx &=&\doint \delta _{a}(x)\overline{%
\delta _{b}(x)}dx \\
&=&\doint \left( \frac{1}{2\beta }\sum_{k\in \Sigma }e^{ika}e^{-ikx}\right)
\left( \frac{1}{2\beta }\sum_{h\in \Sigma }e^{-ihb}e^{ihx}\right) dx \\
&=&\frac{1}{\left( 2\beta \right) ^{2}}\left( \sum_{k,h\in \Sigma
}e^{ika}e^{-ihb}\doint e^{ihx}e^{-ikx}dx\right)  \\
&=&\frac{1}{2\beta }\sum_{k\in \Sigma }e^{ika}e^{-ihb}\delta _{h}^{k}=\frac{1%
}{2\beta }\sum_{k\in \Sigma _{n}}e^{ika}e^{-ikb} \\
&=&\frac{1}{2\beta }\sum_{k\in \Sigma }e^{ik\left( a-b\right) }=\ \underset{%
n\uparrow \sigma }{\lim }\ \frac{1}{2\beta _{n}}\sum_{k\in \Sigma
_{n}}e^{ik\left( b_{n}-a_{n}\right) }
\end{eqnarray*}%
By our definitions we have that%
\begin{eqnarray*}
k &=&l\eta _{n},\ \ l=-n^{2},....,n^{2}-2,n^{2}-1 \\
b-a &=&\ \underset{n\uparrow \sigma }{\lim }\ b_{n}-a_{n} \\
b_{n}-a_{n} &=&L_{n}\eta _{n},\ \ \ L_{n}\in \mathbb{Z};
\end{eqnarray*}%
thus%
\begin{equation*}
\int_{-\beta _{n}}^{\beta _{n}}\delta _{n,a_{n}}(x)\delta _{n,b_{n}}(x)dx=%
\frac{1}{2\beta _{n}}\sum_{k\in \Sigma _{n}}e^{ikL_{n}\eta _{n}}=\frac{1}{%
2\beta _{n}}\sum_{l=-n^{2}}^{n^{2}-1}e^{ilL_{n}\eta _{n}^{2}}.
\end{equation*}%
So, if $b_{n}-a_{n}=L_{n}\eta _{n}=0,$ then $e^{\pi ilL_{n}\eta _{n}^{2}}=1$
and by (\ref{gina})%
\begin{equation*}
\int \delta _{a_{n}}(x)\delta _{b_{n}}(x)dx=\frac{2n^{2}}{2\beta _{n}}=\frac{%
1}{\eta _{n}};
\end{equation*}%
if $b_{n}\neq a_{n},$ since $\eta _{n}^{2}=\pi n^{-2}$ and $L_{n}\in \mathbb{%
Z}\backslash \left\{ 0\right\} ,$ we have that

\begin{eqnarray*}
\int \delta _{a_{n}}(x)\delta _{b_{n}}(x)dx &=&\frac{1}{2\beta _{n}}%
\sum_{l=-n^{2}}^{n^{2}-1}e^{ilL_{n}\eta _{n}^{2}}=\frac{1}{2\beta _{n}}%
\sum_{l=-n^{2}}^{n^{2}-1}e^{\pi ilL_{n}/n^{2}} \\
&=&\frac{e^{i\pi L_{n}}-e^{-i\pi L_{n}}}{e^{\pi iL_{n}/n^{2}}-1} \\
&=&\frac{2i}{e^{\pi iL_{n}/n^{2}}-1}\cdot \sin \left( \pi L_{n}\right) =0
\end{eqnarray*}

Concluding%
\begin{equation*}
\int \delta _{a_{n}}(x)\delta _{b_{n}}(x)dx=\frac{\delta _{a_{n}}^{b_{n}}}{%
\eta _{n}}
\end{equation*}%
and hence%
\begin{equation*}
\doint \delta _{a}(x)\delta _{b}(x)dx=\frac{\delta _{b}^{a}}{\eta }.
\end{equation*}

$\square $

By the above result and the definition of $\delta _{q}(x),$ it turns out that%
\begin{equation*}
\forall q,x\in \Sigma ,\ \delta _{q}(x)=\doint \delta _{q}(t)\delta
_{x}(t)dt=\frac{\delta _{x}^{q}}{\eta }.
\end{equation*}%
Notice the different notation of the Dirac Delta ultrafunction $\delta
_{q}(x)$ and the Kronecker symbol $\delta _{x}^{q}.$

Moreover, we have that 
\begin{equation*}
\left\{ \sqrt{\eta }\delta _{q}(x)\right\} _{q\in \Sigma }
\end{equation*}%
is a orhonormal basis and we will refer to it as to the \textbf{canonical
basis}. Then, using the canonical basis, an ultrafunction $u$ can be
represented as follows:%
\begin{equation}
u(x)=\eta \sum_{q\in \Sigma }\left( \doint u(t)\delta _{q}(t)dt\right)
\delta _{q}(x).  \label{mina}
\end{equation}

So, for every $u\in V_{\sigma }$, it makes sense to write 
\begin{equation}
u(x)=\eta \sum_{q\in \Sigma }u(q)\delta _{q}(x).  \label{cina}
\end{equation}%
and we get the following results:

\begin{corollary}
\label{anna}If $u\in V_{\sigma }$ and $\forall q\in \Sigma ,$ $u(q)=0,$ then 
$\forall x\in \overline{\mathbb{R}}$, $u(x)=0.$
\end{corollary}

\begin{corollary}
\label{annaN}Two tempered ultrafunctions which coincide on $\Sigma $ are
equal.
\end{corollary}

\begin{corollary}
\label{lucia}If $u\in V_{\sigma }$, then%
\begin{equation*}
\doint u(x)dx=\eta \sum_{q\in \Sigma }u(q);
\end{equation*}%
moreover, if $u,v\in V_{\sigma }$, then%
\begin{equation*}
\doint u(x)v(x)dx=\eta \sum_{q\in \Sigma }u(q)v(q);
\end{equation*}
\end{corollary}

Proof: The first inequality follows from (\ref{cina}) and the linearity of $%
\doint $; also, by (\ref{cina}) and Th. \ref{deltadelta}, we have that%
\begin{eqnarray*}
\doint u(x)v(x)dx &=&\doint \left( \eta \sum_{q\in \Sigma }u(q)\delta
_{q}(x)\right) \left( \eta \sum_{r\in \Sigma }v(r)\delta _{r}(x)\right) dx \\
&=&\eta ^{2}\sum_{q,r\in \Sigma }u(q)v(r)\doint \delta _{q}(x)\delta
_{r}(x)dx \\
&=&\eta \sum_{q,r\in \Sigma }u(q)v(r)\delta _{q}^{r}=\eta \sum_{q\in \Sigma
}u(q)v(q).
\end{eqnarray*}

$\square $

So, in the world of tempered ultrafunctions, the integral is equal to the
Riemann sums. However, this fact holds only for the ultrafunctions $u\in
V_{\sigma }$ and not for all the $\sigma $-limit functions.

\bigskip

\section{ The Fourier transform}

\bigskip

\subsection{Definition and main properties}

\bigskip

\begin{definition}
\label{TF}The Fourier transform of an ultrafunction $u\in V_{\sigma },\
\forall k\in \Sigma $, is given by%
\begin{equation*}
\mathfrak{F}\left[ u\right] (k)=\hat{u}(k)=\frac{1}{\sqrt{2\pi }}\doint
u(x)\ e^{-ikx}dx.
\end{equation*}
\end{definition}

\bigskip

By (\ref{cina}) and Corollary \ref{annaN} it is possible to extend $\hat{u}%
(k)$ to an ultrafunction defined on all $\overline{\mathbb{R}}$ setting $%
\forall k\in \overline{\mathbb{R}}$%
\begin{equation*}
\hat{u}(k)=\eta \sum_{h\in \Sigma }\hat{u}(h)\delta _{q}(k)
\end{equation*}

By Corollary \ref{lucia}, $\hat{u}(k),$ $k\in \Sigma $, can be written as
follows:%
\begin{equation}
\hat{u}(k)=\frac{\eta }{\sqrt{2\pi }}\cdot \sum_{x\in \Sigma }u(x)\ e^{-ikx}
\label{pipa}
\end{equation}%
Moreover since $\left\{ \frac{e^{ikx}}{\sqrt{2\beta }}\right\} _{k\in \Sigma
}=\left\{ \sqrt{\frac{\eta }{2\pi }}e^{ikx}\right\} _{k\in \Sigma }$ is an
orthonormal basis, $\hat{u}(k)$ can be regarded, up to the constant $\sqrt{%
\eta }$, as the Fourier component of $u(x)$ in this basis; namely we have
that

\begin{eqnarray*}
u(x) &=&\sum_{k\in \Sigma }\left( \doint u(y)\sqrt{\frac{\eta }{2\pi }}\
e^{-iky}dy\right) \sqrt{\frac{\eta }{2\pi }}e^{ikx} \\
&=&\frac{\eta }{2\pi }\sum_{k\in \Sigma }\left( \doint u(y)\
e^{-iky}dy\right) e^{ikx} \\
&=&\frac{\eta ^{2}}{2\pi }\cdot \sum_{k\in \Sigma }\left( \sum_{y\in \Sigma
}u(y)\ e^{-iky}\right) e^{ikx} \\
&=&\frac{\eta }{\sqrt{2\pi }}\sum_{k\in \Sigma }\hat{u}(k)e^{ikx}
\end{eqnarray*}

\textbf{Examples: }By Th. \ref{deltadelta}, we have that%
\begin{equation}
\doint e^{iqx}e^{-ixk}dx=2\pi \delta _{q}(k).  \label{ganza}
\end{equation}%
and hence%
\begin{equation*}
\mathfrak{F}\left[ e^{iq(\cdot )}\right] (k)=\sqrt{2\pi }\delta _{q}(k)
\end{equation*}%
Moreover, it is immediate to check that%
\begin{equation*}
\mathfrak{F}\left[ \delta _{q}\right] (k)=\frac{1}{\sqrt{2\pi }}e^{-ikq}.
\end{equation*}

Thus the Fourier transform is a change of basis which sends the
trigonometric basis into the canonical basis: 
\begin{equation}
\left\{ \frac{e^{ikx}}{\sqrt{2\beta }}\right\} _{k\in \Sigma }\overset{%
\mathfrak{F}}{\longrightarrow }\left\{ \sqrt{\eta }\delta _{q}(x)\right\}
_{q\in \Sigma }  \label{carla}
\end{equation}%
and hence it is a unitary operator and we have that%
\begin{equation*}
\mathfrak{F}\left[ \delta _{q}\right] (k)=\frac{1}{\sqrt{2\pi }}e^{-ikq}.
\end{equation*}%
Hence we have the following result:

\begin{theorem}
The Fourier transform%
\begin{equation*}
\mathfrak{F}:V_{\sigma }\rightarrow V_{\sigma }
\end{equation*}%
is invertible and we have that 
\begin{equation*}
\mathfrak{F}^{-1}\left[ v\right] (x)=\frac{1}{\sqrt{2\pi }}\doint v(k)\
e^{ixk}dk=\frac{\eta }{\sqrt{2\pi }}\sum_{k\in \Sigma }v(k)\ e^{ikx}.
\end{equation*}%
Moreover%
\begin{equation*}
\doint u(x)\overline{v(x)}dx=\doint \hat{u}(k)\overline{\hat{v}(k)}dk
\end{equation*}
\end{theorem}

\textbf{Proof:} It is an immediate consequence of (\ref{carla}).

$\square $

\bigskip

\subsection{The position operator}

Now let us consider the operator%
\begin{equation*}
\check{x}:V_{\sigma }\rightarrow V_{\sigma }
\end{equation*}%
defined as follows%
\begin{equation*}
\check{x}u(x)=\eta \sum_{q\in \Sigma }qu(q)\delta _{q}(x)
\end{equation*}%
Thus if $x\in \Sigma ,$ $\check{x}u(x)=xu(x);$ however if $x\notin \Sigma ,\ 
$it might happen that $\check{x}u(x)\neq xu(x).$ This possibility is
excluded by (\ref{cina}) if $xu(x)\in V_{\sigma }$, but this fact is false
for some $u\in V_{\sigma }.$ Borrowing this name from quantum mechanics, we
will call $\check{x}$ the \textit{position operator.}

We can extend to the world of ultrafunctions the well known relations
between the position operator and $-iD_{x}$:

\begin{theorem}
\label{pp+}For every $u,v\in V_{\sigma },$

\begin{enumerate}
\item \label{p2+}$\mathfrak{F}\left[ D_{x}u\right] (k)=i\check{k}\widehat{u}%
(k);\ $

\item \label{p2++}$\mathfrak{F}\left[ \check{x}u\right] (k)=iD_{k}\widehat{u}%
(k).$

\item \label{p2}the operators $\check{x}$ and $-iD_{x}$ are Hermitian
operators (in the sense that they are $\sigma $-limit of Hermitian operators
in finite dimensional spaces) and the canonical and the trigonometric bases
are the corresponding bases of eigenvalues.
\end{enumerate}
\end{theorem}

\textbf{Proof}: \ref{p2+} - If $k\in \Sigma $, we have%
\begin{eqnarray*}
\mathfrak{F}\left[ D_{x}u\right] (k) &=&\frac{1}{\sqrt{2\pi }}\doint
D_{x}u(x)\ e^{-ikx}dx=-\frac{1}{\sqrt{2\pi }}\doint u(x)\ D_{x}e^{-ikx}dx \\
&=&\frac{ik}{\sqrt{2\pi }}\doint u(x)\ e^{-ikx}dx=ik\widehat{u}(k)
\end{eqnarray*}%
Thus, for a generic $k\in \overline{\mathbb{R}}$, we have%
\begin{equation*}
\mathfrak{F}\left[ D_{x}u\right] (k)=i\eta \sum_{q\in \Sigma }q\hat{u}%
(q)\delta _{q}(k)=i\check{k}\widehat{u}(k)
\end{equation*}

\ref{p2++} - By Th. \ref{lucia} and (\ref{pipa}), we have that 
\begin{eqnarray*}
\mathfrak{F}\left[ \check{x}u\right] (k) &=&\frac{1}{\sqrt{2\pi }}\doint 
\check{x}u(x)\ e^{-ikx}dx=\frac{\eta }{\sqrt{2\pi }}\sum_{q\in \Sigma
}qu(q)e^{-ikq} \\
&=&\frac{\eta }{\sqrt{2\pi }}iD_{k}\sum_{q\in \Sigma }u(q)e^{-ikq}=iD_{k}%
\widehat{u}(k)
\end{eqnarray*}

\ref{p2} - Trivial.

$\square $

\bigskip

\bigskip

\section{The operator ($%
{{}^\circ}%
$)}

\bigskip

We want to identify ultrafunctions with suitable functions or distributions;
namely, we want to define a liner operator 
\begin{equation*}
(%
{{}^\circ}%
):\mathcal{S}^{\prime }\left( \mathbb{R}\right) \rightarrow V_{\sigma }
\end{equation*}%
such that%
\begin{equation}
\forall \varphi \in \mathcal{S}\left( \mathbb{R}\right) ,\ \forall x\in 
\mathbb{R},\ \varphi 
{{}^\circ}%
(x)=\varphi (x)  \label{nina}
\end{equation}%
\begin{equation*}
\forall T\in \mathcal{S}^{\prime }\left( \mathbb{R}\right) ,\ \forall
\varphi \in \mathcal{S}\left( \mathbb{R}\right) ,\ \doint T%
{{}^\circ}%
(x)\varphi 
{{}^\circ}%
(x)dx=\left\langle T,\varphi \right\rangle
\end{equation*}%
\begin{equation*}
\forall T\in \mathcal{S}^{\prime }\left( \mathbb{R}\right) ,\ DT%
{{}^\circ}%
=\left( DT\right) 
{{}^\circ}%
\end{equation*}

\bigskip

We will realize this program by two steps: first we define $(%
{{}^\circ}%
)$ for $C_{\tau }^{0}\left( \mathbb{R}\right) ,$ secondly we extend it to $%
\mathcal{S}^{\prime }\left( \mathbb{R}\right) ,$ the family of tempered
distributions.

\subsection{The ($%
{{}^\circ}%
$)-operator for continuous slowly increasing functions}

\bigskip

For what we have seen until now, the expression $\doint f(x)\delta _{q}(x)dx$
makes sense when $f$ is an ultrafunction; now we will define it for any $%
f\in C_{\tau }^{0}.$ As it is easy to imagine, we set:%
\begin{equation}
\doint f(x)\delta _{q}(x)dx:=\ \underset{n\uparrow \sigma }{\lim }%
\int_{-\beta _{n}}^{\beta _{n}}f(x)\delta _{q,n}(x)dx  \label{claretta}
\end{equation}%
where $\delta _{q,n}(x)$ is defined by (\ref{pirlo}),

Since $f(y)$ is a slowly increasing function, the above sequence is slowly
increasing and hence the integral converges.

\begin{definition}
\label{luisa}If $f\in C_{\tau }^{0},$ we set%
\begin{equation*}
f%
{{}^\circ}%
(x)=\eta \sum_{q\in \Sigma }\left( \doint f(y)\delta _{q}(y)d\ y\right)
\delta _{q}(x);
\end{equation*}%
where $\ x=\ \underset{n\uparrow \sigma }{\lim }\ x_{n}\in \overline{\mathbb{%
R}}$; in particular, if $x\in \Sigma $,%
\begin{equation}
f%
{{}^\circ}%
(x)=\doint f(y)\delta _{x}(y)d\ y  \label{pipina}
\end{equation}
\end{definition}

\bigskip

The operator $(%
{{}^\circ}%
)$ can also be characterized in the following way:

\bigskip

\begin{theorem}
\label{lalla}The operator:%
\begin{equation*}
(%
{{}^\circ}%
):C_{\tau }^{0}\rightarrow V_{\sigma }
\end{equation*}%
can be written as follows: 
\begin{eqnarray}
f%
{{}^\circ}%
(x) &=&\frac{1}{2\beta }\sum_{k\in \Sigma }\left( \doint f(y)\
e^{-iky}dy\right) e^{ikx}  \label{un} \\
&=&\ \underset{n\uparrow \sigma }{\lim }\ \frac{1}{2\beta _{n}}\sum_{k\in
\Sigma _{n}}\left( \int_{-\beta _{n}}^{\beta _{n}}f(y)\ e^{-iky}dy\right)
e^{ikx_{n}}  \label{du} \\
&=&\frac{1}{2\pi }\doint \left( \doint f(y)\ e^{-iky}dy\right) e^{ikx}dk
\label{tr} \\
&=&\ \underset{n\uparrow \sigma }{\lim }\ \frac{1}{2\pi }\int_{-\beta
_{n}}^{\beta _{n}}\left( \int_{-\beta _{n}}^{\beta _{n}}f(y)\
e^{-iky}dy\right) e^{ikx_{n}}dk  \label{qu}
\end{eqnarray}%
where%
\begin{equation*}
\ x=\ \underset{n\uparrow \sigma }{\lim }\ x_{n}
\end{equation*}
\end{theorem}

\textbf{Proof}: By Def. \ref{luisa}, we have that%
\begin{eqnarray*}
f%
{{}^\circ}%
(x) &=&\eta \sum_{q\in \Sigma }\left( \doint f(y)\delta _{q}(y)d\ y\right)
\delta _{q}(x) \\
&=&\ \underset{n\uparrow \sigma }{\lim }\ \eta _{n}\sum_{q\in \Sigma
_{n}}\left( \int_{-\beta _{n}}^{\beta _{n}}f(y)\delta _{q,n}(y)d\ y\right)
\delta _{q,n}(x) \\
&=&\ \underset{n\uparrow \sigma }{\lim }\ P_{n}(f)
\end{eqnarray*}%
where 
\begin{equation*}
P_{n}(f):=\ \eta _{n}\sum_{q\in \Sigma _{n}}\left( \int_{-\beta _{n}}^{\beta
_{n}}f(y)\delta _{q,n}(y)dy\right) \delta _{q,n}(x)
\end{equation*}%
is the orthogonal projection of $f|_{\left[ \beta _{n},-\beta _{n}\right] }$
(extended by periodicity) on $V_{n}$. Then, representing this projection in
the trigonometric basis we have that 
\begin{equation*}
P_{n}(f):=\ \frac{1}{2\beta _{n}}\sum_{k\in \Sigma _{n}}\left( \int_{-\beta
_{n}}^{\beta _{n}}f(y)\ e^{-iky}dy\right) e^{ikx}
\end{equation*}%
Then, 
\begin{eqnarray*}
f%
{{}^\circ}%
(x) &=&\ \underset{n\uparrow \sigma }{\lim }\ P_{n}(f) \\
&=&\ \underset{n\uparrow \sigma }{\lim }\ \frac{1}{2\beta _{n}}\sum_{k\in
\Sigma _{n}}\left( \int_{-\beta _{n}}^{\beta _{n}}f(y)\ e^{-iky}dy\right)
e^{ikx} \\
&=&\frac{1}{2\beta }\sum_{k\in \Sigma }\left( \doint f(y)\ e^{-iky}dy\right)
e^{ikx}
\end{eqnarray*}%
Thus (\ref{un}) and (\ref{du}) hold. In order to prove (\ref{tr}) we recall
Corollary \ref{lucia} and we get%
\begin{eqnarray*}
f%
{{}^\circ}%
(x) &=&\frac{1}{2\beta }\sum_{k\in \Sigma }\left( \doint f(y)\
e^{-iky}dy\right) e^{ikx} \\
&=&\frac{1}{2\beta \eta }\doint \left( \doint f(y)\ e^{-iky}dy\right)
e^{ikx}dk \\
&=&\frac{1}{2\pi }\doint \left( \doint f(y)\ e^{-iky}dy\right) e^{ikx}dk \\
&=&\ \underset{n\uparrow \sigma }{\lim }\ \frac{1}{2\pi }\int_{-\beta
_{n}}^{\beta _{n}}\left( \int_{-\beta _{n}}^{\beta _{n}}f(y)\
e^{-iky}dy\right) e^{ikx_{n}}dk
\end{eqnarray*}

$\square $

\begin{theorem}
\label{barba}If $f\in C^{0}\left( \mathbb{R}\right) $ is a rapidly
decreasing function, then 
\begin{equation*}
\doint f%
{{}^\circ}%
(x)dx=\int f(x)dx
\end{equation*}
\end{theorem}

\textbf{Proof:} Since 
\begin{equation*}
n\mapsto \int f(x)dx-\int_{-\beta _{n}}^{\beta _{n}}f(x)dx
\end{equation*}%
is a rapidly decreasing sequence, from Th. \ref{peso}, we have that%
\begin{equation*}
\int f(x)dx=\ \underset{n\rightarrow \infty }{\lim }\ \int_{-\beta
_{n}}^{\beta _{n}}f(x)dx=\ \underset{n\uparrow \sigma }{\lim }\ \int_{-\beta
_{n}}^{\beta _{n}}f(x)dx=\doint f%
{{}^\circ}%
(x)dx.
\end{equation*}

$\square $

\bigskip

Next, we are going to check that (\ref{nina}) holds.

\begin{theorem}
\label{palla} If $\varphi $ is a function such that $\hat{\varphi}\in L^{1}(%
\mathbb{R}),$ then, $\forall x\in \mathbb{R}$%
\begin{equation*}
\varphi 
{{}^\circ}%
(x)\sim \varphi (x);
\end{equation*}%
moreover, if $\varphi \in \mathcal{S}(\mathbb{R}),$ then, $\forall x\in 
\mathbb{R}$ 
\begin{equation*}
\varphi 
{{}^\circ}%
(x)=\varphi (x).
\end{equation*}
\end{theorem}

\textbf{Proof}: If $\hat{\varphi}\in L^{1}(\mathbb{R}),$then,%
\begin{equation*}
\varphi (x)=\frac{1}{\sqrt{2\pi }}\ \underset{n\rightarrow \infty }{\lim }%
\int_{-\beta _{n}}^{\beta _{n}}\hat{\varphi}(k)e^{ikx_{n}}dk
\end{equation*}%
and hence, from Th. \ref{peso}, we have that 
\begin{equation*}
\varphi (x)\sim \frac{1}{\sqrt{2\pi }}\ \underset{n\uparrow \sigma }{\lim }\
\int_{-\beta _{n}}^{\beta _{n}}\hat{\varphi}(k)e^{ikx_{n}}dk=\varphi 
{{}^\circ}%
(x)
\end{equation*}

If $\varphi \in \mathcal{S}(\mathbb{R}),\ $we have that 
\begin{equation*}
\varphi (x)=\frac{1}{\sqrt{2\pi }}\int \hat{\varphi}(k)e^{ikx}dk=\frac{1}{%
\sqrt{2\pi }}\int_{-\beta _{n}}^{\beta _{n}}\hat{\varphi}%
(k)e^{ikx_{n}}dk+a_{n}
\end{equation*}%
where $a_{n}$ is rapidly decreasing. On the other hand, 
\begin{equation*}
\hat{\varphi}(k)=\frac{1}{\sqrt{2\pi }}\int \varphi (y)e^{-iky}dy=\frac{1}{%
\sqrt{2\pi }}\int_{-\beta _{n}}^{\beta _{n}}\varphi
(y)e^{-iky}dy+b_{n}\left( k\right)
\end{equation*}%
where $\left\Vert b_{n}\right\Vert _{L^{1}}$ is rapidly decreasing. Then%
\begin{eqnarray*}
\varphi (x) &=&\frac{1}{\sqrt{2\pi }}\int_{-\beta _{n}}^{\beta _{n}}\hat{%
\varphi}(k)e^{ikx_{n}}dk+a_{n} \\
&=&\frac{1}{2\pi }\int_{-\beta _{n}}^{\beta _{n}}\left( \int_{-\beta
_{n}}^{\beta _{n}}\varphi (y)e^{-iky}dy\right) e^{ikx_{n}}dk+\frac{1}{\sqrt{%
2\pi }}\int_{-\beta _{n}}^{\beta _{n}}b_{n}\left( k\right) e^{ikx_{n}}dk+%
\frac{1}{\sqrt{2\pi }}a_{n} \\
&=&\frac{1}{2\pi }\int_{-\beta _{n}}^{\beta _{n}}\left( \int_{-\beta
_{n}}^{\beta _{n}}\varphi (y)e^{-iky}dy\right) e^{ikx_{n}}dk+c_{n}
\end{eqnarray*}%
where $c_{n}$ is rapidly decreasing. Then, by Th. \ref{peso} and Th. \ref%
{lalla}, we have that%
\begin{eqnarray*}
\varphi (x) &=&\ \underset{n\rightarrow \infty }{\lim }\ \frac{1}{2\pi }%
\int_{-\beta _{n}}^{\beta _{n}}\left( \int_{-\beta _{n}}^{\beta _{n}}\varphi
(y)e^{-iky}dy\right) e^{ikx_{n}}dk \\
&=&\ \underset{n\uparrow \sigma }{\lim }\ \frac{1}{2\pi }\int_{-\beta
_{n}}^{\beta _{n}}\left( \int_{-\beta _{n}}^{\beta _{n}}\varphi (y)\
e^{-iky}dy\right) e^{ikx_{n}}dk=\varphi 
{{}^\circ}%
(x)
\end{eqnarray*}

$\square $

\subsection{The ($%
{{}^\circ}%
$)-operator for tempered distributions}

\bigskip

In order to extend $(%
{{}^\circ}%
)$ to $\mathcal{S}^{\prime }\left( \mathbb{R}\right) $ we recall the
following result of Schwartz. If $T\in \mathcal{S}^{\prime }\left( \mathbb{R}%
\right) $, then there exists $f\in C_{\tau }^{0}\left( \mathbb{R}\right) $
and $m\in \mathbb{N}$ such that%
\begin{equation*}
T=D^{m}f
\end{equation*}%
This fact suggests immediately the following definition:

\begin{definition}
\label{martina}Given $T=D^{m}f\in \mathcal{S}^{\prime }\left( \mathbb{R}%
\right) ,$ we set 
\begin{equation*}
T%
{{}^\circ}%
=D^{m}f%
{{}^\circ}%
\end{equation*}
\end{definition}

By the above definition, it follows immediately the following result:

\begin{theorem}
We have that%
\begin{equation}
\forall T\in \mathcal{S}^{\prime }\left( \mathbb{R}\right) ,\ \forall
\varphi \in \mathcal{S}\left( \mathbb{R}\right) ,\ \doint T%
{{}^\circ}%
(x)\varphi 
{{}^\circ}%
(x)dx=\left\langle T,\varphi \right\rangle  \label{pinta}
\end{equation}%
and%
\begin{equation}
\forall T\in \mathcal{S}^{\prime }\left( \mathbb{R}\right) ,\ DT%
{{}^\circ}%
=\left( DT\right) 
{{}^\circ}
\label{santamaria}
\end{equation}
\end{theorem}

\textbf{Proof:} First, let us prove (\ref{pinta}):$\ \forall \varphi \in 
\mathcal{S}\left( \mathbb{R}\right) ,$ 
\begin{eqnarray*}
\doint T%
{{}^\circ}%
(x)\varphi 
{{}^\circ}%
(x)\ dx &=&\doint D^{m}f%
{{}^\circ}%
(x)\varphi 
{{}^\circ}%
(x)\ dx=(-1)^{m}\doint f%
{{}^\circ}%
(x)D^{m}\varphi 
{{}^\circ}%
(x)\ dx \\
&=&(-1)^{m}\ \underset{n\uparrow \sigma }{\lim }\ \int_{-\beta _{n}}^{\beta
_{n}}f(x)D^{m}\varphi (x)\ dx.
\end{eqnarray*}%
Since $f(x)D^{m}\varphi (x)$ is a rapidly decreasing function, 
\begin{equation*}
n\mapsto \int f(x)D^{m}\varphi (x)\ dx-\int_{-\beta _{n}}^{\beta
_{n}}f(x)D^{m}\varphi (x)\ dx
\end{equation*}%
is a rapidly decreasing sequence and hence, by Th. \ref{peso}%
\begin{equation*}
\underset{n\uparrow \sigma }{\lim }\ \int_{-\beta _{n}}^{\beta
_{n}}f(x)D^{m}\varphi (x)\ dx=\ \underset{n\rightarrow \infty }{\lim }\
\int_{-\beta _{n}}^{\beta _{n}}f(x)D^{m}\varphi (x)\ dx
\end{equation*}%
and so%
\begin{eqnarray*}
\doint T%
{{}^\circ}%
(x)\varphi 
{{}^\circ}%
(x)\ dx &=&(-1)^{m}\ \underset{n\rightarrow \infty }{\lim }\ \int_{-\beta
_{n}}^{\beta _{n}}f(x)D^{m}\varphi (x)\ dx \\
&=&(-1)^{m}\int f(x)D^{m}\varphi (x)\ dx=\left\langle T,\varphi
\right\rangle .
\end{eqnarray*}

Now let us prove (\ref{santamaria}). If $T=D^{m}f,$ we have that%
\begin{eqnarray*}
(DT)%
{{}^\circ}
&=&\left( DD^{m}f\right) 
{{}^\circ}%
=\left( D^{m+1}f\right) 
{{}^\circ}%
=D^{m+1}f%
{{}^\circ}
\\
&=&D\left( D^{m}f%
{{}^\circ}%
\right) =D\left( D^{m}f\right) 
{{}^\circ}%
=DT%
{{}^\circ}%
\end{eqnarray*}

$\square $

\bigskip

Some of the previous results can be summarized by theorem \ref{one}.

\bigskip 

\textbf{Proof of Th. \ref{one}:} \ref{11} follows from (\ref{pinta}); \ref%
{22} and \ref{33} follow from Th. \ref{palla}; \ref{44} follows from (\ref%
{santamaria}); \ref{88} follows from Th \ref{barba}, \ref{55} follows from
Th \ref{pizzi}, \ref{66} follows from Def. \ref{TF} and (\ref{pinta}).

\bigskip $\square $

\end{document}